\documentclass[preprint,12pt]{elsarticle}

\usepackage{amssymb}
\usepackage{amsthm}

\journal{Elsevier}

\newtheorem{thm}{Theorem}
\newtheorem{conj}{Conjecture}

\newproof{pf}{Proof}
\newdefinition{defi}{Definition}

\begin{document}

\begin{frontmatter}

\title{Some constructive results on Disjoint Golomb Rulers}


\author[asg]{Xiaodong Xu}
\author[ZS]{Baoxin Xiu}
\author[rvt]{Changjun Fan}
\author[els]{Meilian Liang\corref{cor3}}
\ead{meilianliang@gxu.edu.cn}

\address[asg]{Guangxi Academy of Sciences, Nanning 530007, P.R.China.}
\address[ZS]{School of System Science and Engineering, Sun Yat-sen University, \\Guangzhou 510006, P.R. China}
\address[rvt]{The School of Information System and Management, National University of Defense Technology, Changsha 410073, P.R.China }
\address[els]{School of Mathematics and Information Science, Guangxi University,\\ Nanning 530004, P.R.China.}

\cortext[cor3]{Corresponding author}

\begin{abstract}
A set $\{a_i\:|\: 1\leq i \leq k\}$ of non-negative integers is a Golomb ruler if differences $a_i-a_j$, for any $i \neq j$, are all distinct.
All finite Sidon sets are Golomb rulers, and vice versa.
A set of  $I$ disjoint Golomb rulers (DGR) each being a $J$-subset of $\{1,2,\cdots, n\}$ is called an $(I,J,n)$-DGR. Let $H(I, J)$ be the least positive integer $n$ such that there is an $(I,J,n)$-DGR. In this paper, we propose a series of conjectures on the constructions and structures of DGR. The main conjecture states that if $A$ is any set of positive integers such that $|A| = H(I, J)$, then there are $I$ disjoint Golomb rulers, each being a $J$-subset of $A$,
which generalizes the conjecture proposed by Koml{\'o}s, Sulyok and Szemer{\'e}di in 1975 on the special case $I = 1$.
This main conjecture implies some interesting conjectures on disjoint Golomb rulers.
We also prove some constructive results on DGR, which improve or generalize some basic inequalities on DGR proved by Kl{\o}ve.
\end{abstract}


\end{frontmatter}

\section{Introduction}\label{intro}

A $k$-mark  Golomb ruler is a set of $k$ distinct non-negative integers, also called marks, $\{a_i\:|\: 1\leq i \leq k\}$ such that all differences $a_i-a_j, i \neq j$ are distinct. The difference between the maximal and minimal integer is referred to as the length of the Golomb ruler. Golomb rulers give various important applications in engineering, e.g., the radio-frequency allocation for avoiding third-order interference\cite{Babcock, Golomb1977}, the construction of convolutional or LDPC codes\cite{Robinson2, Chen2012, Kim}, the design of recovery schemes for faulty computers\cite{klonowska2005optimal}, the superresolution multipoint ranging\cite{Oshiga}, etc.

Golomb rulers have been studied by some mathematicians and computer scientists. Various algebraic methods have been proposed to construct Golomb rulers \cite{bose1962theorems, drakakis2009symmetry, ET1941,ruzsa1993solving,singer1938theorem}. However, \textit{Optimal Golomb Rulers}, each of which is the shortest Golomb ruler possible for a given number of marks\cite{Polash}, can only be discovered or verified by exhaustive computer search.  Lengths of $k$-mark optimal rulers have been determined only for $k \leq 27$ so far.

A Sidon sequence (or Sidon set), named after the Hungarian mathematician Simon Sidon, 
is a sequence $A = \{a_1, a_2, ...\}$ of positive integers in which all pairwise sums $a_i + a_j$
$(i \le j)$ are different. Sidon introduced the concept in his investigations of Fourier series \cite{Sidon}.
All finite Sidon sets are Golomb rulers, and vice versa.
Sidon sets have been studied independently in combinatorial number theory for decades.
We can find a survey of Sidon sequences and their generalizations in \cite{O'Bryant2004survey}.

There are some new results on problems related to Golomb rulers or Sidon sets. The topic of infinite Sidon sets contained within sparse random sets of integers was explored in \cite{Fabian2019} and \cite{Kohayakawa2018}. Research on bounds for the maximum size of Sidon sets in a union of intervals was conducted in \cite{Riblet}. Furthermore, modular Golomb rulers have garnered attention from several researchers, with some novel findings on modular Golomb rulers and optical orthogonal codes presented in \cite{MGolombrulers}.

The problem of disjoint Golomb rulers (abbreviated as DGR), a generalization of the Golomb ruler problem, was first considered in \cite{wende1981} in mobile radio-frequency allocation for a collection of areas avoiding third-order interference within each area.
In 1990, Kl{\o}ve conducted an extensive study on DGR in \cite{Klove}, presenting a series of constructions and inequalities through combinatorial and algebraic methods.

In this paper, we present a series of conjectures regarding the constructions and structures of DGR. Our main conjecture generalizes the conjecture proposed by Koml{\'o}s, Sulyok and Szemer{\'e}di in 1975\cite{kssa1975} and leads to several intriguing conjectures related to DGR. Furthermore, we prove some constructive results on DGR, which either improve or generalize certain basic inequalities previously proved by Kl{\o}ve in \cite{Klove}.

The remainder of the paper is organized as follows. Section \ref{secmaincon} proposes a series of conjectures on DGR and establishes the connections between these proposed conjectures. Section \ref{scr} presents some constructive results, which can be viewed as the initial step towards the conjectures proposed in Section \ref{secmaincon}. Section \ref{secmoreconj} introduces an additional conjecture on DGR, which further implies a conjecture on optimal Golomb rulers. Finally, Section \ref{conclu} concludes the paper.

\section{Some Conjectures on Disjoint Golomb Rulers}\label{secmaincon}

We assume that a set of  $I$ DGR each being a $J$-subset of $\{1,2,\cdots, n\}$ is an $(I,J,n)$-DGR. $H(I, J)$ is defined to be the least positive integer $n$ such that there is an $(I,J,n)$-DGR. To determine $H(I,J)$ is an extremely challenging task in some cases. Kl{\o}ve\cite{Klove} proposed a number of constructions for DGR and gave a table of exact values and bounds on $H(I,J)$ for $I \leq 11, J \leq 9$, which were improved and extended through computer search by Shearer \cite{Shearer}.

We propose some conjectures on DGR in Subsection \ref{21} and \ref{22},
and then prove the connections between the proposed conjectures in Subsection \ref{23}.

\subsection{The Main Conjecture}\label{21}

The Golomb ruler problem can be generalized to arbitrary $n$ integers (not necessarily the set $\{1, \cdots,n\}$). Let $A = \{a_1, \cdots, a_n\}$ be an arbitrary finite set of positive integers. Koml{\'o}s, Sulyok and  Szemer{\'e}di  conjectured in \cite{kssa1975} that $A$ must contain an $m$-mark Golomb ruler with $m = (1+o(1))n^{1/2}$. They proved that $m > cn^{1/2}$ for a certain positive constant $c$. The constant was improved by Ruzsa in \cite{RuzsaII1995}.

As pointed out in \cite{kssa1975}, one can expect that the case of the first $n$ positive integers is ``the worst case".
Now, we generalize it to the following conjecture on DGR.

\vspace{0.3cm}
\begin{conj} \label{genefirstn}
Suppose that $I$ and $J$ are positive integers. If $A$ is any set of positive integers such that $|A| = H(I, J)$, then there are $I$ disjoint Golomb rulers, each being a $J$-subset of $A$.
\end{conj}
\vspace{0.3cm}

Conjecture \ref{genefirstn} is the main conjecture in this paper.
It is not difficult to see that if $I=1$, then Conjecture \ref{genefirstn} coincides with the idea in \cite{kssa1975}.

\subsection{More Conjectures on Disjoint Golomb Rulers}\label{22}

Once  Conjecture \ref{genefirstn} is proved, we can establish some important results on DGR. 
The ideas of the following conjecture can be applied in searching upper bounds on $H(I,J)$  for small values of $I$ and $J$.

\vspace{0.3cm}
\begin{conj} \label{difference}
If $I$ and $J$ are integers, and $I \geq 1$, $J \geq 3$, then $H(I+1, J) \leq H(I, J) + J$.
\end{conj}
\vspace{0.3cm}

Among $(I,J,n)$-DGRs, we are very interested in regular ones. An $(I,J,n)$-DGR is regular if $n = IJ$. As we know, there are infinite many regular DGRs (see \cite{Klove}).

We may propose the following conjectures.

\vspace{0.3cm}
\begin{conj} \label{genefirstn2}
If $H(I_0, J) = I_0 J$, then $H(I, J) = IJ$ for any integer $I > I_0$.
\end{conj}
\vspace{0.3cm}

\vspace{0.3cm}
\begin{conj} \label{genefirstnb}
For any Golomb ruler $A_1 \subseteq \{1, 2, \cdots, (I+1)J\}$ such that $|A_1| = J$, if $H(I,J) = IJ$ then there exists a  regular $(I+1,J,(I+1)J)$-DGR containing $A_1$.
\end{conj}
\vspace{0.3cm}

If Conjecture \ref{genefirstn} holds, then Conjecture \ref{genefirstnb} holds too.
Furthermore, it is not difficult to prove by Mathematical Induction that
if Conjecture \ref{genefirstnb} holds
then Conjecture \ref{genefirstn2} holds.

\subsection{Theorems on Conjectures}\label{23}

Imitating the definition of $H(I,J)$, we define $Y(I,J)$ to be
the least positive integer $n$
such that any set of positive integers with $n$ elements,
contains $I$ disjoint $J$-mark Golomb rulers.
Therefore $Y(I,J) \ge H(I,J)$ always holds.

Hence Conjecture \ref{genefirstn} states that
$Y(I,J) = H(I,J)$ for any positive integers $I$ and $J$.
We know that $Y(1,J)$ exists for any positive integer $J$ by a result in \cite{kssa1975}.
In \cite{RuzsaII1995}, a general result on a system of equations
$$a_ {j1}x_1 + \cdots + a_ {jk}x_k = 0, \qquad j = 1, \cdots, J,$$
with integral coefficients $a_ {ji}$ was proved, and
$$Y(1, J) \leq H(1, 5J)$$
was a special case of that general result.
Although this is far from Conjecture \ref{genefirstn}, it is the best known related result now.

For any integers $I \geq 1$ and $J \geq 3$, we
can see that $Y(I,J)$ is finite and
$Y(I,J) \le I Y(1,J)$.
In fact, we can prove the following result on $Y(I+1, J)$.
Note that we can not prove related Conjecture \ref{difference} on
$H(I+1, J)$.

\begin{thm}\label{ThmY}
For any integers $I \geq 1$ and $J \geq 3$
$$Y(I+1, J) \leq Y(I, J) + J.$$
\end{thm}
\begin{pf}

Suppose that $Y(I, J) = n$, $A$ is any set of positive integers and  $|A| = n+J$.
Let $A_0$ be any $J$-mark Golomb ruler in $A$.
Therefore $|A \setminus A_0| = n$.
By the definition of $Y(I,J)$, we can see
there are $I$ disjoint $J$-mark Golomb rulers $A_1, \cdots, A_I$ in $A \setminus A_0$.
So $A_0, A_1, \cdots, A_I$ are $I+1$ disjoint $J$-mark Golomb rulers, all of which are subsets of $A$.
Thus $Y(I+1, J) \leq Y(I, J) + J$.
\end{pf}

By Theorem \ref{ThmY} we can see that $Y(I,J) \leq Y(1,J) + (I - 1)J$.
It is important and difficult to improve the best known upper bound on $Y(1,J)$ given in \cite{RuzsaII1995}.

Now let us prove that if Conjecture \ref{genefirstn} holds,
then both Conjecture \ref{difference} and Conjecture \ref{genefirstn2} hold.
\vspace{0.3cm}

\begin{thm}\label{genefirstn1to2NEW}
(1) If Conjecture \ref{genefirstn} holds, then Conjecture \ref{difference} holds.
(2) If Conjecture \ref{difference} holds, then Conjecture \ref{genefirstn2} holds.
(3) If Conjecture \ref{genefirstn} holds, then Conjecture \ref{genefirstn2} holds.
\end{thm}
\vspace{0.3cm}
\begin{pf}
(1) If Conjecture \ref{genefirstn} holds,
then $H(I+1, J) = Y(I+1, J)$ and $H(I, J) = Y(I, J)$.
By Theorem \ref{ThmY},
$H(I+1, J) = Y(I+1, J) \leq Y(I, J) + J = H(I, J) + J$.
Hence Conjecture \ref{difference} holds when Conjecture \ref{genefirstn} holds.

(2) Suppose that Conjecture \ref{difference} holds and $H(I_0, J) = I_0 J$, then
$(I_0+1)J \leq H(I_0+1, J) \leq H(I_0,J)+J = I_0J+ J = (I_0+1)J$. So $H(I_0+1, J) = (I_0+1)J$.
Thus Conjecture \ref{genefirstn2} holds by Mathematical Induction.

(3) Based on (1) and (2), we have that
if Conjecture \ref{genefirstn} holds, then Conjecture \ref{genefirstn2} holds.
\end{pf}
\vspace{0.3cm}

It is interesting that whether the inequality $H(I+1, J) \leq H(I, J) + J$
in Conjecture \ref{difference} strictly holds when $H(I, J) > IJ$.
In (1) of the proof of Theorem \ref{genefirstn1to2NEW},
we have that
$$H(I+1, J) = Y(I+1, J) \leq Y(I, J) + J = H(I, J) + J$$
when Conjecture \ref{genefirstn} holds.
Therefore if $Y(I+1, J) < Y(I, J) +J$ holds when $H(I, J) > IJ$,
then Conjecture \ref{difference} strictly holds when
Conjecture \ref{genefirstn} holds and
$H(I, J) > IJ$.
Thus it is interesting to know if $Y(I+1, J) < Y(I, J) +J$ holds
whenever $H(I, J) > IJ$, which seems difficult.

\section{Some constructive results}
\label{scr}

Throughout this section, suppose that $a, b, J$ are all positive integers, where $J \geq 3$.

It is evident that $H(a+b,J) \leq H(a,J) + H(b,J)$ (see \cite{Klove}). Furthermore, the following theorem can be easily proved.

\begin{thm}
\label{tI0}
If $H(I_0,J) = I_0J$ for $I_1 \le I_0 \le I_2$, then $H(I,J) = IJ$ for $2I_1 \le I \le 2I_2$.
\end{thm}

Later, we will prove a similar result in Theorem \ref{almosttoregular}.
Now, let us prove some new inequalities on DGR in the following theorem.

\begin{thm}
\label{tHa_b}
(1) $H(a+b,J) \leq H(a,J)+H(b,J-1)+b$, where $a \geq b$. In particular, $H(2a-1,J) \leq H(a,J)+H(a-1,J-1)+a-1$;  \\
(2) $H(2a,J) \leq 2H(a,J-1)+2a$.
\end{thm}
\begin{pf}
(1) Let $F_0= \{X_1, \cdots, X_a\}$ be an $(a,J,n)$-DGR, and
$F_1= \{Y_1, \cdots, Y_b\}$ be a $(b,J-1,m)$-DGR, where $n=H(a,J)$
and $m=H(b,J-1)$. For any $i \in \{1, \cdots, a\}$, let $W_i =
\{x+m \;|\; x \in X_i\}$. Let $F_2= \{W_1, \cdots, W_a\}$. For
any $i \in \{1, \cdots, b\}$, let $U_i = Y_i \cup
\{n+m+i\}$. Let $F_3= \{U_1, \cdots, U_b\}$. Therefore
$F_2 \cup F_3$ is an $(a+b,J,H(a,J)+H(b,J-1)+b)$-DGR, and
$H(a+b,J) \leq H(a,J)+H(b,J-1)+b$.
In particular, if $b=a-1$, then $H(2a-1,J) \leq H(a,J) +H(a-1,J-1)+a-1$.

(2) can be similarly proved.
\end{pf}

Instead of $a \geq b$ in (1) of Theorem \ref{tHa_b},
if $H(a,J) \ge H(b,J-1)$, then $H(a+b,J) \leq H(a,J)+H(b,J-1)+b$ holds.
This can be proved by the same construction in the proof of (1) in Theorem \ref{tHa_b}.

If $H(a,J) = aJ$ and $H(b,J) = bJ$, then $H(a+b,J) = (a+b)J$.
In this case the upper bound on $H(a+b,J) $
given by $ H(a,J) + H(b,J)$ equals to the exact value of $H(a+b,J) $, and can not be improved.

Suppose that $H(a,J) \geq aJ+1$.
To improve $H(a+b,J) \leq H(a,J)+H(b,J)$, now let us prove the following constructive results
in Theorem \ref{gapw}.

The inequalities in Theorem \ref{gapw}, make it possible to construct
regular $(I_1 +I_2, J)$-DGR based on known non-regular $(I_1, J)$-DGR and regular $(I_2, J)$-DGR in some cases.

\begin{thm} \label{gapw}
Suppose that $H(a,J) \geq aJ+1$. If $\{A_i\; |\; 1 \leq i \leq a \}$ is an $(a, J)$-DGR in $\{1, \cdots, H(a,J)\}$,
and there is a gap of size $w$ in the union of $\{A_i\; | \;1 \leq i \leq a\}$, 
then we can prove that $H(a+b,J) \leq H(a,J) + H(b,J) - w$ and
$H(2a,J) \leq 2H(a,J) - 2w$.
\end{thm}
\begin{pf}
Let $H(a,J) =n$ and $H(b,J) = m$.
Suppose that $\{B_i\; |\; 1 \leq i \leq b \}$ is a $(b, J)$-DGR in $\{1, \cdots, m\}$.
Suppose a gap of size $w$ in the union of $\{A_i\; | \;1 \leq i \leq a\}$ is $\{t+i \;|\; 1 \leq i \leq w\}$. 
We may suppose that $t \leq n-w-t$ as well.

Note that $w < H(1,J)$ always hold. 
Otherwise, if $w \geq H(1,J)$, then we can construct an $(a, J)$-DGR in $\{1, \cdots, n-1 \}$ by the following way.

Suppose that $n \in A_{i_0}$, and $D_i = A_i$ for any $i \in \{1, \cdots, a\} \setminus \{i_0\}$. 
Let $D_{i_0}$ be any $J$-ruler in $\{t+1, \cdots, t+w\}$. Therefore we have constructed a $(a, J)$-DGR in $\{1, \cdots, n-1 \}$, what is impossible because that $H(a,J) =n$. Hence $w < H(1,J)$.

Now let us give the constructions in two subcases respectively.

(1) If $n-w-t \leq m$, then we can construct an $(a+b, J)$-DGR $\{C_i\; |\; 1 \leq i \leq a+b \}$ in $\{1, \cdots, n+m-w \}$ as follows.

For $j \in \{1, \cdots, t\}$, $j \in C_i$ if and only if $j \in A_i$, where $i \in \{1, \cdots, a\}$;

For $j \in \{t+x \; |\; 1 \leq x \leq m \}$, $j \in C_{a+i}$ if and only if $x \in B_i$, where $i \in \{1, \cdots, b\}$;

For $j \in \{t+m+y \; |\; 1 \leq y \leq n-t-w \}$, $j \in C_i$ if and only if $y+t+w \in A_i$, where $i \in \{1, \cdots, a\}$.

(2) If $n-w-t > m$, then we can construct an $(a+b, J)$-DGR $\{C_i\; |\; 1 \leq i \leq a+b \}$ in $\{1, \cdots, n+m-w)\}$ as follows.

For $j \in \{1, \cdots, n\}$, $j \in C_i$ if and only if $j \in A_i$, where $i \in \{1, \cdots, a\}$;

For $j \in \{n+x \; |\; 1 \leq x \leq m-w \}$, $j \in C_{a+i}$ if and only if $x \in B_i$, where $i \in \{1, \cdots, b\}$;

For $j \in \{t+y \; |\; 1 \leq y \leq w \}$, $j \in C_i$ if and only if $m-w+y \in B_i$, where $i \in \{1, \cdots, b\}$.

It is not difficult to see that in both cases above
we have constructed an $(a+b, J, n+m-w)$-DGR
$\{C_i\; |\; 1 \leq i \leq a+b \}$.
Therefore
$H(a+b,J) \leq H(a,J) + H(b,J) - w$.

To prove $H(2a,J) \leq 2H(a,J) - 2w$, let us construct
a $(2a, J)$-DGR $\{D_i\; |\; 1 \leq i \leq 2a\}$ in $\{1, \cdots, 2n-2w \}$ as follows.

For $j \in \{1, \cdots, t\}$, $j \in D_i$ if and only if $j \in A_i$, where $i \in \{1, \cdots, a\}$;

For $j \in \{t+x \; |\; 1 \leq x \leq n-t-w \}$, $j \in D_{a+i}$ if and only if $n+1-x \in A_i$, where $i \in \{1, \cdots, a\}$;

For $j \in \{n-w+y \; |\; 1 \leq y \leq n-t-w \}$, $j \in D_i$ if and only if $t+w+y \in A_i$, where $i \in \{1, \cdots, a\}$.

For $j \in \{2n-2w-t+z \; |\; 1 \leq z \leq t \}$, $j \in D_{a+i}$ if and only if $n+1-z \in A_i$, where $i \in \{1, \cdots, a\}$.

By this construction we can see that $\{D_i\; |\; 1 \leq i \leq 2a\}$ is a $(2a, J, 2n-2w)$-DGR.
Therefore $H(2a,J) \leq 2H(a,J) - 2w$.
\end{pf}

If $a$ is small, then the largest gap in the union of $\{A_i\; | \;1 \leq i \leq a\}$ can not be small, because that it is no smaller than $\left\lceil{H(a,J) -aJ \over {aJ-1}}\right\rceil$. Hence we may suppose that $w \ge \left\lceil{{H(a,J) -aJ} \over {aJ-1}}\right\rceil$ in Theorem \ref{gapw}.
Therefore by Theorem \ref{gapw} we have the following theorem.

\begin{thm} \label{gapwcoro}
If $H(a,J) \geq aJ+1$, then
$$H(a+b,J) \leq H(a,J) + H(b,J) - \left\lceil{{H(a,J) -aJ} \over {aJ-1}}\right\rceil,$$
and
$$H(2a,J) \leq 2H(a,J) - 2\left\lceil{{H(a,J) -aJ} \over {aJ-1}}\right\rceil.$$
\end{thm}

Based on Theorem \ref{gapw}, we can prove the following theorem that is better than Theorem \ref{tI0}.

\begin{thm} \label{almosttoregular}
If $I_1 < I_2$,
$H(i,J) \le iJ+1$ for even $i \in \{I_1, \cdots, I_2\}$,
and $H(i,J) = iJ$ for odd $i \in \{I_1, \cdots, I_2\}$,
then $H(I,J) = IJ$ for $I \in \{2I_1, \cdots, 2I_2\}$.
\end{thm}
\begin{pf}
For any even $I \in \{2I_1, \cdots, 2I_2\}$, $H(I,J) = IJ$ by Theorem \ref{gapw}.
Suppose that $I \in \{2I_1, \cdots, 2I_2\}$ and $I$ is odd. Let $I = 2k+1$. Therefore $I_1 \le k \le I_2 - 1$.
If $H(k,J) = kJ$ and $H(k+1,J) = (k+1)J$, then $H(I,J) = (2k+1)J = IJ$.
Otherwise by Theorem \ref{gapw} we know that $H(I,J) \le H(k,J) + H(k+1,J) - 1$.
Because that either $k$ or $k+1$ is odd,
$H(k,J) + H(k+1,J) \le kJ + (k+1)J +1$.
Therefore
$H(k,J) + H(k+1,J) - 1 \le kJ + (k+1)J +1 -1 = (2k+1)J = IJ$. Hence $H(I,J) \le IJ$.
Therefore $H(I,J) = IJ$ for $I \in \{2I_1, \cdots, 2I_2\}$ and we finish the proof.
\end{pf}

Now let us prove an inequality between $H(2,J-1)$ and $H(1,J)$.

\begin{thm} \label{H2j-1}
For any integer $J \geq 3$,
$H(2,J-1) \leq H(1,J)+1$.
\end{thm}
\begin{pf}
Suppose that $n = H(1,J)$ and $A$ is a $J$-mark Golomb ruler
in $\{1, \cdots, n\}$. Let $A_1$ be the set $A+1 = \{a+1 \:| \:a \in A\}$.
Therefore $A_1$ is also a $J$-mark Golomb ruler. 

If $A$ and $A_1$ are disjoint, then $H(2, J-1) < H(2,J) \leq H(1,J)+1$ in this subcase.
Therefore $H(2,J-1) < H(1,J)+1$.

Otherwise, if $A \cap A_1$ is not empty, we can prove that there is at most one common element of $A$ and $A_1$. 
In fact, if $a_1, a_2 \in A \cap A_1$ and $a_1 < a_2$, then $a_1-1, a_2-1, a_1, a_2 \in A$ when $a_2 > a_1+1$, and $a_1-1$, $a_1$,$a_1+1 \in A$ when $a_2 = a_1 +1$, both of which contradict to that 
$A$ is a $J$-mark Golomb ruler. 
Suppose that $a$ is the common element of them.
Since that $A \setminus \{a\}$ and $A_1 \setminus \{H(1,J)+1\}$ are two disjoint
$(J-1)$-mark Golomb rulers, we have that $H(2,J-1) \leq H(1,J) +1$. So we finish the proof. 
\end{pf}

\section{More Conjectures on Disjoint Golomb Rulers and Optimal Golomb Rulers} \label{secmoreconj}

Let $G(k)$ be the length of an optimal Golomb ruler with $k$ marks. Singer proved that if $q$ is a power of a prime, then there exist $q + 1$ integers that have distinct differences modulo $q^2 + q + 1$ and thus
form a Golomb ruler\cite{singer1938theorem}, which implies that $G(k) < k^2-k+1$ if $k-1$ is a prime power.
Based on the Singer's construction,
the following results were proved by Kl{\o}ve in \cite{Klove}.

\vspace{0.3cm}
\begin{thm}
\label{tpp}
If $p$ is a prime power, then
$H(p+1, p) = p^2+p$, $H(p, p-1) \leq p^2-2$, $H(p-1, p) \leq p^2-1$.
\end{thm}
\vspace{0.3cm}

Maybe both
$H(I, I-1) = I^2-I$ and $H(I-1, I) = I^2-I$ hold
for any integer $I \geq 4$.
Furthermore, we propose the following conjecture on $H(I, I+2)$.

\vspace{0.3cm}
\begin{conj} \label{HII2regular}
For any integer $I \geq 4$, $H(I, I+2) = I(I+2)$.
\end{conj}
\vspace{0.3cm}

By the results in \cite{Klove} and \cite{Shearer}, we know that Conjecture \ref{HII2regular} holds for $I \in \{4,5,6,7\}$.

If both Conjecture \ref{genefirstn2} and Conjecture \ref{HII2regular} hold, then $H(I, J) = IJ$ for any integer $I \geq J-2$ and $J \geq 6$.
This seems far from reach, and the best known related result is the following theorem of Kl{\o}ve proved in \cite{Klove}.

\vspace{0.3cm}
\begin{thm} \label{Klovetheorem}\cite{Klove}
Denote $\iota(J)$ the minimum integer such that $H(I, J) = IJ$ for all $I \ge \iota(J)$.
For all $J \ge 2$, we have
$${{H(1,J)-2} \over {J-1}} < \iota(J) \le 4p-1$$
where $p$ is any odd prime such that $p \ge J-1$.
\end{thm}
\vspace{0.3cm}

Note that Kl{\o}ve proved the upper bound for $\iota(J)$ in Theorem \ref{Klovetheorem} constructively.

Moreover, if both Conjecture \ref{genefirstn} and Conjecture \ref{HII2regular} hold, for any two disjoint $J$-mark Golomb rulers 
in $\{1, 2,  \cdots, IJ\}$, $g_1$ and $g_2$, there is a regular $(I,J,IJ)$-DGR containing $g_1$ and $g_2$,
where $6 \leq J \leq I$.
It may be time-consuming to computationally verify this even for small $I$ and $J$.

Upper bounds on $G(k)$
given by Singer's work in \cite{singer1938theorem}
are only for the cases in which $p$ is a prime power.
If Conjecture \ref{HII2regular} holds, it is not difficult to see that the following bounds on $G(k)$ hold for general situations.

\vspace{0.3cm}
\begin{conj} \label{SRkk+2b}
$G(k+2) < k^2 + k$ for any $k \geq 6$.
\end{conj}
\vspace{0.3cm}

For $k \ge 6$, Conjecture \ref{SRkk+2b} is stronger than the conjecture that $G(k) < k^2$ for all $k > 0$,
which was first mentioned by Erd\H{o}s in an equivalent form\cite{ET1941} and is recognized to be true for $k \leq 65000$ up to now\cite{dimitromanolakis2002analysis}.

\section{Conclusions and Remarks} \label{conclu}

The problem of finding disjoint Golomb rulers is
an interesting generalization of the Golomb ruler probelm.
We have generalized the DGR problem to arbitrary $n$ positive integers and conjectured
that if $A$ is any set of positive integers such that $|A| = H(I, J)$, then there are $I$ disjoint Golomb rulers, each being a $J$-subset of $A$.
We have proved that
this conjecture implies some important conjectures on DGRs.

Theoretical research and computational verification of the conjs proposed in this paper are our future tasks.
These tasks may be difficult.

\section*{Acknowledgements}
This work is partially supported by the NSF of China (71471176, 11361008).


\end{document}